\documentclass[a4paper,12pt]{amsart}
\usepackage{amsfonts}
\usepackage{amssymb}
\usepackage{ifthen}
\usepackage{graphicx}
\nonstopmode \numberwithin{equation}{section}
\usepackage[usenames]{color}
\newtheorem{thm}{Theorem}[section]
\newtheorem{lem}{Lemma}[section]
\newtheorem{cor}{Corollary}[section]

\newtheorem{cl}{Claim}[section]
\newtheorem{ca}{Case}[section]
\newtheorem{sca}{Subcase}[section]
\newtheorem{scl}[section]{Subclaim}
\newtheorem{conj}[equation]{Conjecture}

\theoremstyle{definition}
\newtheorem{defn}{Definition}[section]
\newtheorem{op}[equation]{Open Problem}
\newtheorem{ques}[equation]{Question}
\newtheorem{rem}{Remark}[section]
\newtheorem{exam}[equation]{Example}

\newcounter {own}
\def\theown {\thesection       .\arabic{own}}

\newenvironment{pf}[1][]{%
 \vskip 3mm
 \noindent
 \ifthenelse{\equal{#1}{}}%
  {{\slshape Proof. }}%
  {{\slshape #1.} }%
 }%
{\qed\bigskip}

\newcounter{alphabet}
\newcounter{tmp}
\newenvironment{Thm}[1][]{\refstepcounter{alphabet}%
\bigskip%
\noindent%
{\bf Theorem \Alph{alphabet}}%
\ifthenelse{\equal{#1}{}}{}{ (#1)}%
{\bf .} \itshape}{\vskip 8pt}

\makeatletter
\newcommand{\Ref}[1]{\@ifundefined{r@#1}{}{\setcounter{tmp}{\ref{#1}}\Alph{tmp}}}
\makeatother

\newenvironment{Lem}[1][]{\refstepcounter{alphabet}%
\bigskip%
\noindent%
{\bf Lemma \Alph{alphabet}}%
{\bf .} \itshape}{\vskip 8pt}

\newcommand{\diam}{{\operatorname{diam}}}

\newcommand{\dist}{{\operatorname{dist}}}

\def\be{\begin{equation}}
\def\ee{\end{equation}}

\newcommand{\ben}{\begin{enumerate}}
\newcommand{\een}{\end{enumerate}}

\newcommand{\blem}{\begin{lem}}
\newcommand{\elem}{\end{lem}}
\newcommand{\bthm}{\begin{thm}}
\newcommand{\ethm}{\end{thm}}
\newcommand{\bcor}{\begin{cor}}
\newcommand{\ecor}{\end{cor}}
\newcommand{\beg}{\begin{exam}}
\newcommand{\eeg}{\end{exam}}
\newcommand{\begs}{\begin{examples}}
\newcommand{\eegs}{\end{examples}}
\newcommand{\bdefe}{\begin{defn}}
\newcommand{\edefe}{\end{defn}}
\newcommand{\bprob}{\begin{prob}}
\newcommand{\eprob}{\end{prob}}
\newcommand{\bques}{\begin{ques}}
\newcommand{\eques}{\end{ques}}
\newcommand{\bei}{\begin{itemize}}
\newcommand{\eei}{\end{itemize}}
\newcommand{\bcon}{\begin{conj}}
\newcommand{\econ}{\end{conj}}
\newcommand{\bop}{\begin{op}}
\newcommand{\eop}{\end{op}}

\newcommand{\bca}{\begin{ca}}
\newcommand{\eca}{\end{ca}}
\newcommand{\bsca}{\begin{sca}}
\newcommand{\esca}{\end{sca}}

\newcommand{\bcl}{\begin{cl}}
\newcommand{\ecl}{\end{cl}}

\newcommand{\bscl}{\begin{scl}}
\newcommand{\escl}{\end{scl}}

\newcommand{\bcons}{\begin{conjs}}
\newcommand{\econs}{\end{conjs}}
\newcommand{\bprop}{\begin{propo}}
\newcommand{\eprop}{\end{propo}}
\newcommand{\br}{\begin{rem}}
\newcommand{\er}{\end{rem}}
\newcommand{\brs}{\begin{rems}}
\newcommand{\ers}{\end{rems}}
\newcommand{\bo}{\begin{obser}}
\newcommand{\eo}{\end{obser}}
\newcommand{\bos}{\begin{obsers}}
\newcommand{\eos}{\end{obsers}}
\newcommand{\bpf}{\begin{pf}}
\newcommand{\epf}{\end{pf}}
\newcommand{\ba}{\begin{array}}
\newcommand{\ea}{\end{array}}
\newcommand{\beq}{\begin{eqnarray}}
\newcommand{\beqq}{\begin{eqnarray*}}
\newcommand{\eeq}{\end{eqnarray}}
\newcommand{\eeqq}{\end{eqnarray*}}

\newcounter{minutes}
\divide\time by 60
\newcounter{hours}
\multiply\time by 60 \addtocounter{minutes}{-\time}

\begin{document}

\bibliographystyle{amsplain}

\title{Quasim\"obius invariance of uniform domains}

\def\thefootnote{}
\footnotetext{ \texttt{\tiny File:~\jobname .tex,
          printed: \number\year-\number\month-\number\day,
          \thehours.\ifnum\theminutes<10{0}\fi\theminutes}
} \makeatletter\def\thefootnote{\@arabic\c@footnote}\makeatother

\author{Qingshan Zhou${}^{\mathbf{*}}$}
\address{Qingshan Zhou, School of Mathematics and Big Data, Foshan university,  Foshan, Guangdong 528000, People's Republic
of China} \email{qszhou1989@163.com; q476308142@qq.com}

\author{Antti Rasila}
\address{Antti Rasila, Technion -- Israel Institute of Technology, Guangdong Technion, 241 Daxue Road, Shantou, Guangdong 515063, People's Republic of China} \email{antti.rasila@gtiit.edu.cn; antti.rasila@iki.fi}

\date{}
 \keywords{Uniform domain, min-max property,
(relative) quasim\"{o}bius mapping, quasisymmetric mapping, natural condition.\\
${}^{\mathbf{*}}$ Corresponding author}

\begin{abstract}
In this paper, we study quasim\"obius invariance of uniform domains in Banach spaces. We first investigate implications of certain geometric properties of domains in Banach spaces, such as the (diameter) uniformity, the $\delta$-uniformity and the min-max property. Then we show that all of these conditions are equivalent if the domain is $\psi$-natural. As applications, we answer partially to an open question proposed by V\"ais\"al\"a, and  provide a new method to prove a recent result in \cite{HLVW}, which also gives an answer to another question raised by V\"ais\"al\"a.
\end{abstract}

\thanks{Qingshan Zhou was supported by NSF of China (No. 11901090), and by Department of Education of Guangdong Province, China (Grant No. 2018KQNCX285). Antti Rasila was supported by NNSF of China (No. 11971124).}

\maketitle{} \pagestyle{myheadings} \markboth{}{Quasim\"obius invariance of uniform domains}

\section{Introduction and main results}\label{sec-1}
Many results of the classical function theory have their counterparts in the setting of quasiconformal mappings in $n$-dimensional Euclidean spaces. To further extend the scope of this theory, V\"ais\"al\"a developed from late 1980's onwards the theory of (dimension) free quasiconformal mappings in Banach spaces \cite{Vai-1, Vai-2, Vai-3, Vai-4, Vai-5}. The main advantage of this approach is in avoiding use of volume integrals and conformal modulus, which allows one to study the quasiconformality of mappings in infinite dimensional Banach spaces and other metric spaces without volume measures. Recently, this line of research has a subject of several investigations (see e.g. \cite{HL,HLVW,HRWZ,LPZ,LVZ17,rt2,WZGHR} and references therein).

In this paper, we consider the relationship between uniform domains and relative quasim\"obius maps in real Banach spaces. The main objective is to study two open questions formulated by V\"ais\"al\"a in \cite{Vai-2} and \cite{Vai-5}. We will give a partial solution to one of these questions and provide a new method to answer the other one. Following \cite{Vai-5}, we assume that $E$ is a real Banach space with dimension at least two, a proper domain $G\subsetneq E$ is a nonempty connected open set, and $d_G(x)=\dist(x,\partial G)$ for $x\in G$. Let us begin with the definition of uniform domains in Banach spaces.

\bdefe\label{def1.1} Let $E$ be a real Banach space and $G \subsetneq E$ a domain, and let $c\geq 1$. We say that $G$ is $c$-{\it uniform} if each pair of points $x$, $y$ in $G$ can be connected by an arc $\gamma$ in $G$ satisfying:
\begin{enumerate}
\item $\ell(\gamma)\leq c\,|x-y|$, and
\item $\min\{\ell(\gamma[x,z]),\ell(\gamma[z,y])\}\leq c\,d_G(z)$ for all $z\in \gamma$,
\end{enumerate}
\noindent where $\ell(\gamma)$ denotes the length of $\gamma$, $\gamma[x,z]$ the part of $\gamma$ between $x$ and $z$.
At this time, $\gamma$ is said to be a $c$-{\it uniform arc}.
\edefe
\br In 1978, Martio and Sarvas \cite{MS} introduced the twisted interior cone condition in connection to showing global injectivity properties for locally injective mappings. Since then, many other characterizations of uniform  domains have been established, see \cite{GH, GO, Martio-80}. Uniform domains can be understood as a class of domains developed in the context of generalizing the Riemann mapping theorem for quasiconformal maps in $\mathbb{R}^n$ with $n\geq 3$, a question that still remains open. This class of domains has numerous geometric and function theoretic properties that make it useful in many fields of modern mathematical analysis as well (see e.g. \cite{GP,Jo80,Vai-0,Vai,Vai-5}).
\er

In \cite{Martio-80}, Martio studied the quasiconformal invariance of uniform domains in $\overline{R}^n=R^n\cup\{\infty\}$. He obtained the following result by showing several equivalent conditions for uniform domains.

\begin{Thm}$($\cite[Theorems 5.4 and 6.2]{Martio-80}$)$ Let $c,K\geq 1$. If $G\subsetneq \overline{R}^n$ is a $c$-uniform domain and $f:\overline{R}^n\to \overline{R}^n$ is $K$-quasiconformal, then $f(G)$ is $c_1$-uniform with $c_1$ depending only on $c,K$ and $n$.
\end{Thm}

Subsequently, V\"ais\"al\"a \cite{Vai-0} introduced the concept of {\it quasim\"obius} maps (see Definition \ref{def2'}) and investigated the relation between this class of mappings and quasiconformal on uniform domains.

\begin{Thm}\label{T-1}$($\cite[Theorem 5.6]{Vai-0}$)$ Let $c,K\geq 1$ and $n\geq 2$. If $G\subsetneq \overline{R}^n$ is a $c$-uniform domain and $f:G\to G'\subsetneq \overline{R}^n$ is a $K$-quasiconformal homeomorphism, then $G'$ is $c_1$-uniform if and only if $f$ is $\eta$-quasim\"obius.
\end{Thm}

Moreover, V\"ais\"al\"a has generalized Theorem \Ref{T-1} to real Banach spaces by using the free quasiconformality (FQC) theory in \cite{Vai-2}.

\begin{Thm}\label{T-2}$($\cite[Theorem 7.18]{Vai-2}$)$ Let $G$ and $G'$ be domains in Banach spaces $E$ and $E'$, respectively. If $G$ is $c$-uniform and $f:G\to G'$ is $\varphi$-FQC, then $G'$ is $c_1$-uniform if and only if $f$ is $\eta$-quasim\"obius.
\end{Thm}

The main motivation of this paper is the following two questions posed by V\"ais\"al\"a, the second of which is still open. Note that the definition of $(M,C)$-CQH and relative quasim\"obius mapping see Definitions \ref{def2'} and \ref{def3'}.

\begin{ques}\label{q-r1}$($\cite[Question 7.19]{Vai-2}$)$ Suppose that $G$ is a $c$-uniform domain and that $f:G\to G'$ is $(M,C)$-CQH. If $f$ extends to a homeomorphism $\overline{f}:\overline{G}\to \overline{G'}$ and $\overline{f}$ is $\theta$-quasim\"obius relative to the boundary $\partial G$, is $G'$ $c'$-uniform?
\end{ques}

\br
It is worth of mentioning that Huang, Li, Vuorinen, and Wang \cite{HLVW} answered Question \ref{q-r1} affirmatively. Their proofs are based on several concepts and results in the free quasiworld \cite{Vai-5}; such as coarse quasihyperbolic length and solid arcs, the equivalence of uniform and $\varphi$-uniform domains, the diameter cigar theorems for uniform domains.
\er

\begin{ques}\label{q-r2}$($\cite[13.2.10]{Vai-5}$)$ Suppose that $G$ is a $c$-uniform domain and that a homeomorphism ${f}:\overline{G}\to \overline{G'}$ is $\eta$-quasim\"obius relative to $\partial G$ and maps $G$ onto $G'$. Is $G'$ $c'$-uniform with $c'=c'(c,\eta)$?
\end{ques}

We focus our attention to Question \ref{q-r2}. Indeed, Question \ref{q-r2} is more difficult than Question \ref{q-r1} because one does not assume that $f$ is a coarsely quasihyperbolic map. To deal with this question, we first prove the following implications between certain geometric properties for domains in Banach spaces, such as the diameter uniformity, the min-max property and the $\delta$-uniformity for some $0<\delta<1$. Also see Subsection \ref{sub-domain} for the definitions.

\begin{thm}\label{thm-main1} Suppose that $G\subsetneq E$ is a domain, then we have the following implications: $G$ is $c$-uniform $\Rightarrow$ $G$ satisfies the min-max property $\Rightarrow$ $G$ is diameter $c_1$-uniform $\Leftrightarrow$ $G$ is $\delta$-uniform for some $0<\delta<1$.
\end{thm}
\br The min-max property for domains in $\mathbb{R}^n$ was introduced by Gehring and Hag in \cite{GH}. They extended the properties of hyperbolic geodesic in $\mathbb{B}^n$ to more general domains, and studied the relationship between this property, the uniformity and quasiconformal extension property. In \cite{Martio-80}, Martio introduced the concept of the $\delta$-uniformity property in terms of the cross-ratio of four points for domains in $\mathbb{R}^n$. By using this condition, he obtained certain general properties of uniform domains.
\er

As an application of Theorem \ref{thm-main1}, we prove the following quasim\"obius invariance of diameter uniform  domains and $\delta$-uniform ($0<\delta<1$) domains.

\bcor\label{c-1} Let $G$ and $G'$ be proper domains in Banach spaces $E$ and $E'$, respectively. Let $f:G\to G'$ be a $\theta$-quasim\"obius homeomorphism.
\begin{enumerate}
  \item If $G$ is diameter $c$-uniform, then $G'$ is diameter $c'$-uniform with $c'=c'(c,\theta)$;
  \item If $G$ is $\delta$-uniform with $0<\delta< 1$, then $G'$ is $\delta'$-uniform with $\delta'=\delta'(\delta,\theta)\in (0,1)$.
\end{enumerate}
\ecor

Next, we consider the converse of Theorem \ref{thm-main1}. One easily see that $c$-uniform domain is diameter $c$-uniform. For the converse, it follows from \cite[Theorem 4.5]{Martio-80} that a diameter $c$-uniform domain in $\mathbb{R}^n$ is $c_1$-uniform with $c_1$ depending on $c$ and $n$. It is natural to ask whether this holds in infinite-dimensional Banach spaces. Our second main result is an attempt in this direction.

\begin{thm}\label{thm-main2} Let $G\subsetneq E$ be a domain. Then $G$ is $c$-uniform if and only if $G$ is diameter $c_1$-uniform and $\psi$-natural.
\end{thm}

\br The definition of $\psi$-natural domains is given in Definition \ref{n-1}.  We note that every proper domain in $\mathbb{R}^n$ is $\psi$-natural with $\psi=\psi(n)$ (see \cite[Corollary 2.18]{Vu1}). In an infinite dimensional Hilbert space, the broken tube construction  in \cite{Vai04} provides an example of a domain, which is not natural. Moreover, it is diameter $c$-uniform but not $c_1$-uniform for any constant $c_1\geq 1$.
\er

There are several applications of Theorems \ref{thm-main1} and \ref{thm-main2} in studying the relationship between uniform domains and (relative) quasim\"obius maps. The following is an immediate one because a $\psi$-uniform domain (see (\ref{z-3}) for the definition) is $\psi$-natural.

\bcor\label{c-2} Let $G\subsetneq E$ be a domain. Then $G$ is $c$-uniform if and only if $G$ is diameter $c_1$-uniform and $\psi$-uniform. In particular, a diameter uniform convex domain in a Banach space is uniform.
\ecor

As the second application, we answer partially to Question \ref{q-r2} as follows.
\begin{thm}\label{thm-main3} Let $G$ and $G'$ be proper domains in Banach spaces $E$ and $E'$, respectively. Suppose that $G$ is $c$-uniform and that a homeomorphism ${f}:\overline{G}\to \overline{G'}$ is $\eta$-quasim\"obius relative to $\partial G$ and maps $G$ onto $G'$, then $G'$ is diameter $c'$-uniform with $c'=c'(c,\eta)$. If in addition $G'$ is $\psi$-natural, then $G'$ is $c_1'$-uniform with $c'_1=c'_1(c,\eta,\psi)$.
\end{thm}

Moreover, we show that Theorem \ref{thm-main3} gives a way to affirmatively answer Question \ref{q-r1}. Our method is quite different from the one used in \cite{HLVW}.

\begin{thm}\label{thm-main4} Let $G$ and $G'$ be proper domains in Banach spaces $E$ and $E'$, respectively. Suppose that $f:G\to G'$ is $(M,C)$-CQH, and that $f$ extends to a homeomorphism ${f}:\overline{G}\to \overline{G'}$ which is $\eta$-quasim\"obius relative to  $\partial G$.
\begin{enumerate}
  \item If $G$ is $\psi$-natural, then $G'$ is $\psi'$-natural with $\psi'=\psi'(\psi,M,C,\eta)$,
  \item If $G$ is $c$-uniform, then $G'$ is $c'$-uniform with $c'=c'(c,M,C,\eta)$.
\end{enumerate}
\end{thm}

By Theorems \ref{thm-main2} and \ref{thm-main4}, we further obtain the following consequence.

\bcor\label{c-3} Let $G$ and $G'$ be proper domains in Banach spaces $E$ and $E'$, respectively. Suppose that $f:G\to G'$ is $(M,C)$-CQH, and that $f$ extends to a homeomorphism ${f}:\overline{G}\to \overline{G'}$ which is $\eta$-quasim\"obius relative to  $\partial G$. If $G$ is diameter $c$-uniform and $G'$ is $\psi$-natural, then both $G$ and $G'$ are $c_1$-uniform with $c_1=c_1(c,M,C,\theta,\psi)\geq 1$.
\ecor

The rest of this paper is organized as follows. In Section \ref{sec-2}, we recall necessary definitions and preliminary results. The proofs of Theorems \ref{thm-main1} and \ref{thm-main2} are given in Section \ref{sec-3}. Section \ref{sec-4} is devoted to the proofs of Theorem \ref{thm-main3} and \ref{thm-main4}, and the proofs of Corollaries \ref{c-1}, \ref{c-2} and \ref{c-3} are presented in Section \ref{sec-5}.

\section{Preliminaries and notations}\label{sec-2}

\subsection{\bf Notation.}

Let letters $A,B,C,...$ denote positive numerical constants. Similarly, $C(a,b,c,...)$ denotes universal positive functions of the parameters $a,b,c,...$. Sometimes we write $C=C(a,b,c,...)$ to emphasize the parameters on which $C$ depends and abbreviate $C(a,b,c,...)$ to $C$.

Following the notation and terminology of \cite{HLVW,Vai-5}, we use $E$ and $E'$ to denote real Banach spaces with dimension at least $2$. The norm of a vector $x$ in $E$ is written as $|x|$, and for every pair of points $z_1$, $z_2$ in $E$, the distance between them is denoted by $|z_1-z_2|$, the closed line segment with endpoints $z_1$ and $z_2$ by $[z_1, z_2]$. The one-point extension of $E$ is the Hausdorff space $\dot{E}=E\cup\{\infty\}$, where the neighborhoods of $\infty$ are the complements of closed bounded sets of $E$.

For a set $A$ in $E$, we use $\overline{A}$ to denote the completion of $A$ and $\partial A=\overline{A}\setminus A$ to be its norm boundary. For a bounded set $A$ in $E$, $\diam \,A$ is the diameter of $A$. Let
$$B(x,r)=\{ z\in E:\; |z-x|<r\},\; \overline{B}(x,r)=\{ z\in E:\; |z-x|\leq r\},$$
and $S(x,r)=\{ z\in E:\; |z-x|=r\}$.

Let $X$ be a metric space. A curve is a continuous function $\gamma:$ $[a,b]\to X$. The length of $\gamma$ is defined by
$$\ell(\gamma)=\sup\Big\{\sum_{i=1}^{n}|\gamma(t_i)-\gamma(t_{i-1})|\Big\},$$
where the supremum is taken over all partitions $a=t_0<t_1<t_2<\ldots<t_n=b$. The curve is called {\it rectifiable} if $\ell(\gamma)<\infty$. The metric space $X$ is called {\it rectifiably connected} if each pair of points can be connected by a rectifiable curve.

The length function associated with a rectifiable curve $\gamma$: $[a,b]\to X$ is $s_{\gamma}$: $[a,b]\to [0, \ell(\gamma)]$, defined by
$s_{\gamma}(t)=\ell(\gamma|_{[a,t]})$ for $t\in [a,b]$. For any rectifiable curve $\gamma:$ $[a,b]\to X$, there is a unique parametrization $\gamma_s:$ $[0, \ell(\gamma)]\to X$ such that $\gamma=\gamma_s\circ s_{\gamma}$. Obviously, $\ell(\gamma_s|_{[0,t]})=t$ for $t\in [0, \ell(\gamma)]$. The parametrization $\gamma_s$ is called the {\it arclength parametrization} of $\gamma$. For a rectifiable curve $\gamma$ in $X$, the line integral over $\gamma$ of each Borel function $\varrho:$ $X\to [0, \infty)$ is
$$\int_{\gamma}\varrho ds=\int_{0}^{\ell(\gamma)}\varrho\circ \gamma_s(t) dt.$$

\subsection{\bf  Domains in Banach spaces.}\label{sub-domain}

In this part, we assume that $E$ is a real Banach space with dimension at least two, and $G\subsetneq E$ is a domain. We begin with the definition of the quasihyperbolic metric. Note that this metric was first introduced by Gehring and Palka \cite{GP} in the case of proper domains in $R^n$. It has been recently used by many authors in the study of quasiconformal mappings and related questions, see e.g. \cite{GO,LVZ2,LVZ3,rt2,ZRL}.

Recall that the {\it quasihyperbolic length} of a rectifiable curve $\alpha$ in the norm metric in $G$ is the number (cf. \cite{GO,GP}):
$$\ell_k(\alpha)=\int_{\alpha}\frac{|dz|}{d_{G}(z)}.
$$

For each pair of points $z_1$, $z_2$ in $G$, the {\it quasihyperbolic distance}
$k_G(z_1,z_2)$ between $z_1$ and $z_2$ is defined in the usual way:
$$k_G(z_1,z_2)=\inf\ell_k(\alpha),
$$
where the infimum is taken over all rectifiable curves $\alpha$ joining $z_1$ and $z_2$ in $G$. An arc $\alpha$ from $z_1$ to $z_2$ is a {\it quasihyperbolic geodesic} if
$$\ell_k(\alpha)=k_G(z_1,z_2).$$
It is known that a quasihyperbolic geodesic between every pair of points in $G$ exists if the dimension of $E$ is finite, see \cite[Lemma 1]{GO}. However, this is not true in infinite dimensional Banach spaces (cf. \cite[Example 2.9]{Vai-1}).

Let us remark that the second author with Talponen further investigated in \cite{rt2} properties of quasihyperbolic geodesics in Banach spaces. They demonstrated that in a strictly convex Banach space with the Radon-Nikodym property, the quasihyperbolic geodesics are unique. Moreover, they provided an example to show that for a convex domain in a non-reflexive Banach space, it is possible that there is no quasihyperbolic geodesic between any given pair of points in the domain.

In order to overcome this shortcoming, V\"ais\"al\"a introduced the following concept in \cite{Vai-2}.

\bdefe Let $G\subsetneq E$ be a domain and let $c\geq 1$. An arc $\alpha\subseteq G$ is a $c$-{\it neargeodesic} if $\ell_k(\alpha[x,y])\leq c\;k_G(x,y)$ for all $x, y\in \alpha$.
\edefe

Moreover, V\"ais\"al\"a \cite{Vai-2} proved the following property concerning the existence of neargeodesics in the domains of Banach spaces.

\begin{Thm}\label{LemA} $($\cite[Theorem 3.3]{Vai-2}$)$ Let $G\subsetneq E$ be a domain. Then for all points $\{z_1,\, z_2\}\subseteq G$ and for any $c>1$, there is a $c$-neargeodesic $\alpha$ joining $z_1$ and $z_2$ in $G$.
\end{Thm}

We also record the following elementary inequality (cf. \cite{Vai-5}):
\be\label{eq(0000)} k_{G}(z_1, z_2)\geq \Big|\log \frac{d_{G}(z_2)}{d_{G}(z_1)}\Big|,\ee
for all $z_1$, $z_2$ in $G$.

We recall the definition of $\varphi$-uniform domains presented in \cite{Vu2}.

\bdefe
Let $\varphi:[0,\infty)\to [0,\infty)$ be a strictly increasing homeomorphism. A domain $G\subsetneq E$ is called $\varphi$-{\it uniform} if for all $x$, $y$ in $G$, we have
\be\label{z-3}  k_G(x,y)\leq \varphi(r_G(x,y)),\;\;\;\;\;\mbox{where}\;\;\;\;\;r_G(x,y)=\frac{|x-y|}{\min\{d_G(x),d_G(y)\}}.\ee
\edefe

In order to give a simple criterion for $\varphi$-uniform domains,
consider domains  $G$ satisfying the following property \cite[Examples~2.50~(1)]{Vu2}:
there exists a constant $C\ge 1$ such that each pair of points
$x,y\in G$ can be joined by a rectifiable curve $\gamma\in G$ with
$\ell(\gamma)\le C\,|x-y|$ and $\min\{d_G(x),d_G(y)\}\le
C\,\dist(\gamma,\partial G)$. Then $G$ is $\varphi$-uniform with
$\varphi(t)=C^2t$. In particular, every convex domain is
$\varphi$-uniform with $\varphi(t)=t$. However, in general, (unbounded) convex
domains need not be uniform.

Suppose that $\emptyset\not=A\subseteq G\subsetneq E$. We write $$r_G(A)=\sup\{r_G(x,y): x\in A, y\in A\}$$
and
 $$k_G(A)=\sup\{k_G(x,y): x\in A, y\in A\}.$$

\bdefe\label{n-1}  Let $\psi:[0,\infty)\to [0,\infty)$ be an increasing function. A domain $G \subsetneq E$ is called $\psi$-{\it natural} if
$$k_G(A)\leq \psi(r_G(A))$$
for every nonempty connected set $A\subseteq G$ with $r_G(A)<\infty$.
\edefe

This definition is from the paper of V\"ais\"al\"a \cite{Vai-3}, and also studied by Vuorinen in \cite{Vu1}. We note that a $\varphi$-uniform domain is $\varphi$-natural, and every convex domain is $\psi$-natural with $\psi(t)=t$ (see, \cite[Theorems 2.8 and 2.9]{Vai-3}).
In fact, the next two results show that the class of natural domains is fairly large.

\begin{Lem}\label{lem-2}$($\cite[Corollary 2.18]{Vu1}$)$
Every proper domain in $\mathbb{R}^n$ is $\psi_n$-natural with $\psi_n$ depending only on $n$.
\end{Lem}

\begin{Lem}\label{lem-3}$($\cite[Theorem 2.8]{Vai-3}$)$ Let $c\geq 1$. A $c$-uniform domain $G$ in a Banach space $E$ is $\psi$-natural with $\psi=\psi(c)$.
\end{Lem}

But in an infinite dimensional Hilbert space, the broken tube construction in \cite[2.3]{Vai04} provides an example of a domain, which is not natural.

\bdefe Let $G \subsetneq E$ be a domain and let $c\geq 1$. We say that $G$ is $c$-{\it John} if each pair of points $x$, $y$ in $G$ can be connected by an arc $\gamma$ in $G$ such that
$$\min\{\ell(\gamma[x,z]),\ell(\gamma[z,y])\}\leq c\,d_{G}(z),$$
for all $z\in \gamma$. Moreover, the arc $\gamma$ is called a {\it $c$-cone arc}.
\edefe
\br The concept of John domains in Euclidean spaces was first introduced in 1961 by John \cite{Jo61} in connection with his work in elasticity. Recently, Li, Vuorinen and the first author studied several equivalent conditions for John metric spaces in \cite{LVZ17}. Indeed, their method is used in our proof of Theorem \ref{thm-main2}.
\er

\bdefe Let $G \subsetneq E$ be a domain and let $c\geq 1$. $G$ is called {\it diameter $c$-uniform}, if each pair of points $x_{1},x_{2}$ in $G$ can be joined by an arc $\alpha$ in $G$ satisfying:
\begin{enumerate}
\item $\min_{j=1,2}\;\diam\; (\alpha [x_j, x])\leq c\,d_{G}(x)
$ for all $x\in \alpha$, and
\item $\diam \,\alpha \leq c\,|x_{1}-x_{2}|$.
\end{enumerate}
Moreover, $\alpha$ satisfying the above conditions is said to be a {\it diameter uniform arc}.
\edefe
\br It is easy to see that a $c$-uniform domain is diameter $c$-uniform. For the converse, it follows from \cite[Theorem 4.5]{Martio-80} that a diameter $c$-uniform domain in $\mathbb{R}^n$ is $c_1$-uniform with $c_1$ depending on $c$ and $n$. Note that in Banach spaces, V\"ais\"al\"a \cite{Vai04} constructed a broken tube domain which is diameter $c$-uniform but not $c_1$-uniform for any constant $c_1\geq 1$.
\er

\bdefe Let $G \subsetneq E$ be a domain. We say that $G$ has the {\it min-max property} if there exists a family of curves $\Gamma$ in $G$ and a constant $c\geq 1$ such that any pair of points in $G$ can be joined by a curve $\gamma\in \Gamma$ and so that
\be\label{r-1} \frac{1}{c}\min_{j=1,2} |x_j-y| \leq |x-y| \leq c\max_{j=1,2} |x_j-y|, \ee
for each ordered triple of points $x_1,x,x_2\in \gamma$ and each $y\in \partial G$.
\edefe

\bdefe Let $G \subsetneq E$ be a domain and let $0<\delta<1$. $G$ is called $\delta$-{\it uniform} if each pair of points $x_{1},x_{2}$ in $G$ can be joined by an arc $\alpha$ in $G$ such that the cross ratio
\be\label{r-2} \tau(x,x_i,y,x_j)=\frac{|x-y|}{|x-x_i|} \cdot \frac{|x_i-x_j|}{|x_j-y|}\geq \delta,\;\;\;i\neq j\in\{1,2\}, \ee
for all $x\in \alpha\setminus\{x_1,x_2\}$ and $y\in E\setminus G$.
\edefe

\subsection{\bf  Mappings on metric spaces.}\label{sub-1}
 Let $X$ be a metric space and $\dot{X}=X\cup \{\infty\}$. By a
triple in $X$ we mean an ordered sequence $T=(x,y,z)$ of three
distinct points in $X$. The ratio of $T$ is the number
$$\rho(T)=\frac{|y-x|}{|z-x|}.$$ If $f: X\to Y$ is  an injective
map, the image of a triple  $T=(x,y,z)$  is the triple
$f(T)=(f(x),f(y),f(z))$.

Suppose that  $A\subseteq X$. A triple  $T=(x,y,z)$ in $X$ is said to
be a triple in the pair $(X, A)$ if $x\in A$ or if $\{y,z\}\subseteq
A$. Equivalently, both $|y-x|$ and $|z-x|$ are distances from a
point in $A$.

\bdefe \label{def1-0} Let $X$ and $Y$ be two metric spaces, and let
$\eta: [0, \infty)\to [0, \infty)$ be a homeomorphism. Suppose
$A\subset X$. An embedding $f: X\to Y$ is said to be {\it
$\eta$-quasisymmetric} relative to $A$, or $\eta$-$QS$ rel
$A$, if $\rho(f(T))\leq \eta(\rho(T))$  for each triple $T$ in
$(X,A)$. \edefe

Notice that an embedding $f: X\to Y$ is $\eta$-$QS$ relative to  $A$ if
and only if $\rho(T)\leq t$ implies that $\rho(f(T))\leq \eta(t)$
for each triple $T$ in $(X,A)$ and $t\geq 0$ (cf. \cite{Vai-5}).
Obviously, quasisymmetric relative to  $X$ is equivalent to usual quasisymmetric.

A quadruple in $X$ is an ordered sequence $Q=(x,y,z,w)$ of four
distinct points in $X$. The cross ratio of $Q$ is defined to be the
number
$$\tau(Q)=\tau(x,y,z,w)=\frac{|x-z|}{|x-y|}\cdot\frac{|y-w|}{|z-w|}.$$ Observe that the definition can be extended in
the usual manner to the case where one of the points is
$\infty$. For example,
$$\tau(x,y,z,\infty)= \frac{|x-z|}{|x-y|}.$$
If $X_0 \subseteq \dot{X}$ and if $f: X_0\to \dot{Y}$
is an injective map, the image of a quadruple $Q$ in $X_0$ is the
quadruple $f(Q)=(f(x),f(y),f(z),f(w))$. Suppose that $A\subseteq X_0$. We say
that a quadruple $Q=(x,y,z,w)$ in $X_0$ is a quadruple in the pair
$(X_0, A)$ if $\{x,w\}\subseteq A$ or $\{y,z\}\subseteq A$.
Equivalently, all four distances in the definition of $\tau(Q)$ are
(at least formally) distances from a point in $A$.

\bdefe \label{def2'} Let ${X}$ and ${Y}$ be two metric
spaces and let $\eta: [0, \infty)\to [0, \infty)$ be a
homeomorphism. Suppose $A\subseteq \dot{X}$. An embedding $f:
\dot{X}\to \dot{Y}$ is said to be {\it $\eta$-quasim\"obius}
relative to $A$, or $\eta$-$QM$ rel $A$, if the inequality
$\tau(f(Q))\leq \eta(\tau(Q))$ holds for each quadruple in $(X,A)$.
\edefe

Apparently, $\eta$-$QM$ relative to  $X$ is equivalent to $\eta$-quasim\"obius.

We conclude this section by recalling the definition of CQH maps and by presenting a useful result.

\bdefe\label{def3'} Let $G$ and $G'$ be domains in Banach spaces $E$ and $E'$, respectively. Let $M\geq 1$ and $C\geq 0$. A homeomorphism $f:G\to G'$ is said to be $C$-coarsely $M$-quasihyperbolic, or $(M,C)$-CQH if for all $x,y\in G$, we have
$$\frac{k_G(x,y)-C}{M}\leq k_{G'}(f(x),f(y))\leq Mk_G(x,y)+C.$$
\edefe

\begin{Lem}\label{z-2}$($\cite[Lemma 2.14]{Vai-3}$)$  Suppose that ${f}:\overline{G}\to \overline{G'}$ is $\eta$-quasim\"obius relative to  $\partial G$ and that $f(G)=G'$. Suppose also that $x,y\in G$ with $r_G(x,y)\leq t$. Then $r_{G'}(f(x),f(y))\leq \mu(t,\eta)<\infty$.
\end{Lem}

\section{Proofs of Theorems \ref{thm-main1} and \ref{thm-main2}}\label{sec-3}
\subsection{} In this subsection, we prove Theorem \ref{thm-main1} by dividing the arguments into several lemmas.
We use the methods of the proof from \cite[Theorem 2.7]{GH} and \cite[Theorem 5.4]{Martio-80}.

\begin{lem}\label{r-1.0}
Suppose that $G\subsetneq E$ is a $c$-uniform domain, then $G$ satisfies the min-max property.
\end{lem}
\bpf Set
$$\Gamma:=\{\alpha_{u,v}\;\;|\;\;  u,v\in G,\mbox{$\alpha_{u,v}$ is a $2$-neargeodesic with end points $u,v$ in $G$} \}.$$
It follows from Theorem \Ref{LemA} that $\Gamma$ is not empty. Fix $\alpha\in \Gamma$ and for an ordered triple of points $x_1,x,x_2\in \gamma$ and $y\in \partial G$, it suffices to verify (\ref{r-1}) for these points.

First, by \cite[Theorem 10.17]{Vai-5}, we know that $G$ is quasihyperbolic $c_1$-uniform with $c_1=c_1(c)$; for the definition see \cite[10.2]{Vai-5}. Moreover, it follows from Cigar Theorem \cite[10.9]{Vai-5} that the subarc $\gamma[x_1,x_2]$ is $c_2$-uniform with $c_2=c_2(c)$. Thus we have
$$\min_{j=1,2}|x_j-x|\leq \min_{j=1,2} \ell(\gamma[x_j,x])\leq c_2 d_G(x)\leq c_2|x-y|.$$
This yields
\be\label{r-1.1}\min_{j=1,2}|x_j-y|\leq |x-y|+\min_{j=1,2} |x_j-x|\leq (1+c_2)|x-y|.\ee
On the other hand, by the uniformity of $\gamma[x_1,x_2]$, we obtain
\be\label{r-1.2} |x-y|\leq |x-x_1|+|x_1-y|\leq c_2|x_1-x_2|+|x_1-y|\leq (2c_2+1)\max_{j=1,2}|x_j-y|.\ee

Therefore, we see from (\ref{r-1.1}) and (\ref{r-1.2}) that $G$ satisfies the min-max property for this family of curves $\Gamma$. Hence we prove Lemma \ref{r-1.0}.
\epf

\begin{lem}\label{r-2.0}
Suppose that a domain $G\subsetneq E$ satisfies the min-max property, then $G$ is diameter $c_1$-uniform.
\end{lem}
\bpf Fix $x_1,x_2\in G$. If $|x_1-x_2|\leq \frac{1}{2}d_G(x_1)$, in this case,  we see that the line segment $[x_1,x_2]$ is the desired diameter uniform arc.

If $|x_1-x_2|> \frac{1}{2}d_G(x_1)$, take $y\in \partial G$ with $|x_1-y|\leq 2d_G(x_1)$. Then we have
\be\label{r-2.1} |x_1-y|\leq 4|x_1-x_2|.\ee

Because $G$ satisfies the min-max property, we know that there exists a curve $\gamma$ joining $x_1$ and $x_2$ in $G$ and a constant $c\geq 1$ such that (\ref{r-1}) holds. Then for all $x\in \gamma$, we obtain by (\ref{r-2.1}) that
$$|x-y|\leq c\max_{j=1,2} |x_j-y|\leq c(|x_1-x_2|+|x_1-y|)\leq 5c|x_1-x_2|,$$
which implies
\be\label{r-2.2}\diam \;\gamma\leq 10c |x_1-x_2|.\ee

Fix $z\in \gamma$ and take another point $z\in \partial G$ with $|x-z|\leq 2d_G(x)$. Moreover, we claim that for some $j=1$ or $j=2$,
\be\label{r-2.3} \gamma[x_j,x]\subseteq B(z,2cd_G(x)).\ee
Otherwise, there are two points $u_j\in \gamma[x_j,x]$ for $j=1,2$ such that
$$\min_{j=1,2}|u_j-z|>2c d_G(x) \geq c |x-z|,$$
which contradicts (\ref{r-1}). Therefore, we see from (\ref{r-2.3}) that
\be\label{r-2.4} \min_{j=1,2} \diam \;\gamma[x_j,x]\leq 4cd_G(x).\ee
Hence it follows from (\ref{r-2.2}) and (\ref{r-2.4}) that $\gamma$ is the required in this case,  which completes the proof of Lemma \ref{r-2.0}.
\epf

\begin{lem}\label{r-3.0} Let $G\subsetneq E$ be a diameter $c$-uniform domain with $c\geq 1$. Then $G$ is $\delta$-uniform for some constant $\delta=\delta(c)\in(0,1)$.
\end{lem}
\bpf Fix two distinct points $x_1,x_2\in G$. Because $G$ is diameter $c$-uniform, we can choose a diameter $c$-uniform arc $\gamma$ connecting $x_1$ and $x_2$ in $G$. Thus for all $x\in \gamma$ and $y\in E\setminus G$, we may assume without loss of generality that
$$\diam\;\gamma[x_1,x]\leq \diam\;\gamma[x_2,x].$$
In order to show that $G$ is $\delta$-uniform, we need to find lower bounds for the cross ratios $\tau(x,x_1,y,x_2)$ and $\tau(x,x_2,y,x_1)$.

We first consider two cases to check that
\be\label{r-3.0.0} \tau(x,x_1,y,x_2)\geq \frac{1}{2c^2}.\ee
If $|x_2-y|\leq 2c|x_1-x_2|$, we have
$$\tau(x,x_1,y,x_2)=\frac{|x-y||x_1-x_2|}{|x-x_1||x_2-y|} \geq \frac{\diam\;\gamma[x_1,x] \;|x_1-x_2|}{c|x-x_1||x_2-y|}\geq \frac{1}{2c^2}.$$
If $|x_2-y|> 2c|x_1-x_2|$, we observe from the diameter uniformity of $\gamma$ that
$$\max_{j=1,2}|x_j-x|\leq \diam \;\gamma[x_1,x_2] \leq c|x_1-x_2|,$$
and thus
$$\frac{|x-y|}{|x_2-y|}\geq \frac{|x_2-y|-|x-x_2|}{|x_2-y|}\geq \frac{1}{2}.$$
This implies that
$$\tau(x,x_1,y,x_2)=\frac{|x-y||x_1-x_2|}{|x-x_1||x_2-y|}\geq \frac{1}{2c},$$
and we obtain (\ref{r-3.0.0}).

It remains to find a lower bound for $\tau(x,x_2,y,x_1)$. Towards this end, we check that
\be\label{r-3.1} |x-y|\geq \frac{1}{c+1}|y-x_1|.\ee
Suppose on the contrary that $|x-y|< \frac{1}{c+1}|y-x_1|$. Then we have
$$|x_1-x|\geq |y-x_1|-|y-x|>\frac{c}{c+1}|y-x_1| $$
and thus
$$|x-y|>\frac{1}{c} \diam\gamma[x_1,x]\geq \frac{1}{c}|x-x_1|> \frac{1}{c+1}|y-x_1|,$$
which is a contradiction, showing (\ref{r-3.1}). Moreover, we obtain from (\ref{r-3.1}) that
\be\label{r-3.1.1} \tau(x,x_2,y,x_1)=\frac{|x-y||x_1-x_2|}{|x-x_2||x_1-y|} \geq \frac{1}{c(c+1)}.\ee

Therefore, by (\ref{r-3.0.0}) and (\ref{r-3.1.1}) we find that $G$ is $\delta$-uniform with $\delta=1/(2c^2+2c)$. Hence the proof of Lemma \ref{r-3.0} is complete.
\epf

\begin{lem}\label{r-3.0.1} Let $G\subsetneq E$ be a $\delta$-uniform domain with $\delta\in(0,1)$. Then $G$ is diameter $c$-uniform for some constant $c=c(\delta)\geq 1$.
\end{lem}
\bpf Fix two distinct points $u,v\in G$ and take $t=\frac{1}{16}|u-v|$. If $|u-v|<\max\{d_G(u),d_G(v)\}$, then we see that the line segment $[u,v]\subseteq G$ is a diameter uniform arc. Thus we may assume that
\be\label{r-3.00} |u-v|\geq\max\{d_G(u),d_G(v)\}.\ee
Because  $G$ is $\delta$-uniform, we can join $u$ and $v$ by an arc $\alpha$ satisfying the $\delta$-uniformity condition (\ref{r-2}). Denote by $\alpha_0$ the subarc of $\alpha$ containing in $G\setminus \big(B(u,t)\cup B(v,t)\big)$ and  such that $\alpha_0$ connects the spheres $S(u,t)$ and $S(v,t)$. Let $u_0=\alpha_0\cap S(u,t)$ and $v_0=\alpha_0\cap S(v,t)$.

Next, we are going to construct a curve $\beta$ joining $u$ to $u_0$ and satisfying
\be\label{r-3.2} \beta[u,x]\subseteq B(x,C d_G(x)) \ee
for all $x\in \beta$, where $C=C(\delta)\geq 1$ will be determined below.

To this end, choose a sequence of points $u_i$, such that $u_i\in \alpha \cap S(u,t/{2^i})$ for all $i\geq 1$. Moreover, for each $i\geq 0$, pick some point $y_i\in \partial G$ with
$$|y_i-u_i|\leq 2d_G(u_i).$$
We consider two possibilities.

If there is some index $i\geq 0$ such that $|v-y_i|<\frac{1}{2}|u-v|$, thus we have
\beq\nonumber  d_G(u)&\geq& d_G(u_i)-|u-u_i|
\\ \nonumber&\geq&  \frac{1}{2} |u_i-y_i|- |u-u_i|
\\ \nonumber&\geq&  \frac{1}{2}|u-v|-\frac{3}{2}|u-u_i|-\frac{1}{2}|v-y_i|
\\ \nonumber&>& \frac{1}{4}|u-v|-\frac{3}{2}|u-u_0|
\\ \nonumber&>& 2|u-u_0|,
\eeq
where the last equality follows from the choice of $u_0$ in $\alpha$. From the above it follows that we can take the line segment $[u,u_0]=\beta$ joining $u$ to $u_0$ such that (\ref{r-3.2}) holds for $C=1$, as desired.

We are thus left to assume that $|v-y_i|\geq \frac{1}{2}|u-v|$ for all $i\geq 0$. In this case for each $i\geq 0$, we take an arc $\beta_i$ connecting $u_i$ to $u_{i+1}$ satisfying the $\delta$-uniformity condition. Thus $\beta:=\bigcup_{i=0} \beta_i$ is a curve joining $u$ to $u_0$. In order to show that $\beta$ is the required curve, we need some estimates.

By our  choices of the points $u_i$ and $y_i$, the $\delta$-uniformity of $\alpha$ implies that
\be\label{r-3.3} d_G(u_i)\geq \frac{1}{2}|y_i-u_i|\geq \frac{\delta |v-y_i|}{2|u-v|}|u_i-u|=\frac{\delta}{2^{i+5}}|v-y_i| \geq \frac{\delta}{2^{i+6}}|u-v|=:r_i,\ee
for all $i\geq 0$.

Moreover, for all $x\in \beta_i$ and  for every $y\in \partial G$, we claim that
\be\label{r-3.4} |x-y|\geq \frac{\delta^3}{2^{i+12}}|u-v|.\ee
This can be seen as follows. If $x$ lies in the union of $B(u_i,\frac{1}{2}r_i)$ and $B(u_{i+1},\frac{1}{2}r_{i+1})$, by (\ref{r-3.3}) we have
$$|x-y|\geq |y-u_{i+1}|-|x-u_{i+1}|\geq \frac{1}{2}r_{i+1}=\frac{\delta}{2^{i+8}}|u-v|,$$
as needed.

If $x$ lies outside of $B(u_i,\frac{1}{2}r_i)$ and $B(u_{i+1},\frac{1}{2}r_{i+1})$, we have $|x-u_i|\geq \frac{1}{2}r_i$, $|y-u_{i+1}|\geq r_{i+1}$ and
$$|u_i-u_{i+1}|\leq |u_i-u|+|u_{i+1}-u|\leq \frac{t}{2^{i-1}}.$$
Because $\beta_i$ satisfies the $\delta$-uniformity property, the estimates from above produce
$$ |x-y|\geq \frac{\delta |x-u_i||y-u_{i+1}|}{|u_i-u_{i+1}|} \geq \frac{\delta^3}{2^{i+12}}|u-v|,$$
as desired. This shows (\ref{r-3.4}).

We verify (\ref{r-3.2}) by using (\ref{r-3.4}). Fix $x\in \beta$. Then there is an index $k=0,1,\ldots$ such that $x\in \beta_k$. For all $z\in \beta_i$ with $i\geq k$, the $\delta$-uniformity of $\beta_i$ yields
$$|z-u_i|\leq \frac{1}{\delta}|u_i-u_{i+1}|\leq \frac{t}{\delta 2^{i-1}},$$
and therefore
$$|z-u|\leq |z-u_i|+|u_i-u|\leq \frac{t}{\delta 2^{i-2}}.$$
Consequently, by (\ref{r-3.4}) we find that
$$|z-x|\leq |z-u|+|u-x|\leq \frac{t}{\delta 2^{k-3}}\leq \frac{2^{11}}{\delta^4}d_G(x),$$
which implies (\ref{r-3.2}) by letting $C=2^{11}\delta^{-4}$.

A similar arguments as above shows that there is a curve $\gamma$ connecting $v_0$ and $v$ such that
\be\label{r-3.5} \gamma[v,x]\subset B(x,C d_G(x)) \ee
for all $x\in \gamma$.
 
Let $\alpha'=\beta\cup \alpha_0\cup \gamma$. We shall prove that $\alpha'$ satisfies the desired diameter uniform condition. For all $x\in \beta$, we obtain from (\ref{r-3.00}) and (\ref{r-3.2}) that
$$\frac{1}{C}|x-u_0|\leq d_G(u_0) \leq d_G(u)+|u-u_0|\leq 2|u-v|,$$
and thus
\be\label{r-3.6} \diam\; \beta \leq 4C|u-v|.\ee
Similarly, it follows from (\ref{r-3.00}) and (\ref{r-3.5}) that
\be\label{r-3.7} \diam\; \gamma \leq 4C|u-v|.\ee
Moreover, for all $x\in \alpha_0$, by the $\delta$-uniformity of $\alpha$, we have $|u-x|\leq |u-v|/\delta$, which implies
\be\label{r-3.8} \diam\; \alpha_0 \leq \frac{2}{\delta}|u-v|.\ee
Therefore, we see from (\ref{r-3.6}), (\ref{r-3.7}) and (\ref{r-3.8}) that
\be\label{r-3.9}\diam\; \alpha'\leq \diam\;\beta+\diam\;\alpha_0+\diam\;\gamma\leq (8C+\frac{2}{\delta})|u-v|.\ee

It follows from (\ref{r-3.2}), (\ref{r-3.5}) and (\ref{r-3.9}) that we only need to show that there exists a constant $A\geq 1$ such that
$$\min\{\diam \;\alpha'[u,x],\diam \;\alpha'[v,x]\} \leq A d_G(x),$$
for all $x\in\alpha_0$. Because $\alpha$ satisfies the $\delta$-uniformity condition, for all $y\in \partial G$ we thus have
\beq\nonumber |x-y|&\geq& \max\Big\{ \frac{\delta |x-u|}{|u-v|}|y-v|,\frac{\delta |x-v|}{|u-v|}|y-u|\Big\}
\\ \nonumber &\geq&  \frac{\delta}{16}\max\{|y-v|,|y-u|\}
\\ \nonumber &\geq&  \frac{\delta}{32}|u-v|.
\eeq
This together with (\ref{r-3.9}) yields
$$\diam\;\alpha'\leq \Big(8C+\frac{2}{\delta}\Big)\frac{32}{\delta}d_G(x),$$
as desired. Hence we obtain  Lemma \ref{r-3.0.1}.
\epf

\emph{Proof of Theorem \ref{thm-main1}.} From Lemmas \ref{r-1.0}, \ref{r-2.0}, \ref{r-3.0} and \ref{r-3.0.1}, we immediately obtain Theorem \ref{thm-main1}. \qed

\subsection{} In this part, we give the proof of Theorem \ref{thm-main2}. Note that the necessity part follows from Lemma \Ref{lem-3}. To prove the sufficiency, we first develop certain intermediate results. The idea of our arguments for Lemmas \ref{s-1} and \ref{s-2} below comes from \cite{LVZ17}.

\begin{lem}\label{s-1} Assume that $G\varsubsetneq E$ is a $\psi$-natural and diameter $c$-uniform domain. If $x,x_0\in G$ with $d_G(x_0)\leq 2d_G(x)$ and $|x-x_0|\leq c_0d_G(x_0)$ for some constant $c_0\geq 1$, then there is an arc $\alpha$ in $G$ connecting $x$ and $x_0$ such that
$$\ell_k(\alpha[x,x_0])\leq 2\psi(2cc_0(1+c)),$$
where $k$ is the quasihyperbolic metric of $G$.
\end{lem}
\bpf
Because $G$ is diameter $c$-uniform, there is a diameter $c$-uniform arc $\gamma$ in $G$ joining $x$ and $x_0$. Because $d_G(x_0)\leq 2d_G(x)$, we thus have
\beq\label{s-1.1} \frac{1}{2}d_G(x_0) &\leq& \min\{d_G(x), d_G(x_0)\}
\\ \nonumber &\leq&  d_G(z)+\min\{|z-x|, |z-x_0|\}
\\ \nonumber &\leq&  (1+c)d_G(z),
\eeq
for all $z\in\gamma$. Moreover, because $G$ is $\psi$-natural and $|x-x_0|\leq c_0d_G(x_0)$, we obtain by (\ref{s-1.1}) that
\beq\label{z-1}
k(x,x_0) &\leq&  \diam_k(\gamma)\leq \psi\Big(\frac{\diam\;\gamma}{\dist(\gamma,\partial G)}\Big)
\\ \nonumber &\leq& \psi\Big(\frac{2c(1+c)|x-x_0|}{d_G(x_0)}\Big)
\\ \nonumber &\leq& \psi\Big(2cc_0(1+c)\Big),
\eeq
where $\diam_k(\gamma)$ denotes the quasihyperbolic diameter of  $\gamma$.

Next, choose an arc $\alpha$ connecting $x$ and $x_0$ in $G$ such that $\ell_k(\alpha)\leq 2k(x,x_0)$. Therefore, (\ref{z-1}) yields
$$\ell_k(\alpha)\leq 2k(x,x_0)\leq 2\psi(2cc_0(1+c)),$$
as desired. This shows Lemma \ref{s-1}.
\epf

\begin{lem}\label{s-2} Assume that $G\varsubsetneq E$ is a $\psi$-natural and diameter $c$-uniform domain. If $x,x_0\in G$ with $d_G(x_0)> 2d_G(x)$, then there is an arc $\alpha$ in $G$ connecting $x$ and $x_0$ such that
$$\ell_k(\alpha[x,y])\leq 2\psi(8c^2(1+c)),$$
where $y$ is the first point of $\alpha$ such that $d_G(y)=2d_G(x)$ when traversing $\alpha$ from $x$ to $x_0$ and $k$ is the quasihyperbolic metric of $G$.
\end{lem}
\bpf Because $G$ is diameter $c$-uniform, it follows that there is a diameter $c$-uniform arc $\gamma$ in $G$ joining $x$ and $x_0$. Denote by $y_0$ the first point of $\gamma$ such that
$$d_G(y_0)=2d_G(x)<d_G(x_0).$$
Let $\beta=\gamma[x,y_0]$ and $\beta'=\gamma[y_0,x_0]$. Then we consider two possibilities and show that
\be\label{s-2.1} k(x,y_0)\leq \psi(8c^2(1+c)).\ee

If $\diam\;\beta\leq \diam\;\beta'$, thus by the uniformity of $\gamma$ we have
$$\diam\;\beta\leq cd_G(y_0)=2c d_G(x)$$
and
$$d_G(x)\leq |x-z|+d_G(z)\leq (c+1)d_G(z),$$
for all $z\in\beta$. Because $G$ is $\psi$-natural, this implies that
$$k(x,y_0)\leq \diam_k(\beta)\leq \psi\Big(\frac{\diam\;\beta}{\dist(\beta,\partial G)}\Big)\leq \psi(2c(1+c)),$$
as required.

If $\diam\;\beta> \diam\;\beta'$, choose $u\in\beta$ with $\diam\;\beta[x,u]=\diam\;\beta/2$. Again we see from the uniformity of $\beta$ that
$$|x-x_0|\leq 2\diam\;\beta\leq 4\diam\;\beta[x,u]\leq 4cd_G(u)\leq 8cd_G(x),$$
and
\beq \nonumber d_G(x) &=&  \min\{d_G(x),d_G(x_0)\}
\\ \nonumber &\leq&  d_G(z)+\min\{|z-x|, |z-x_0|\}
\\ \nonumber &\leq& (c+1)d_G(z),
\eeq
for any $z\in\beta$. Because $G$ is $\psi$-natural, this yields
\beq \nonumber k(x,y_0) &\leq&  \diam_k(\beta)\leq \psi\Big(\frac{\diam\;\beta}{\dist(\beta,\partial G)}\Big)
\\ \nonumber &\leq&  \psi\Big(\frac{c(c+1)|x-x_0|}{d_G(x)}\Big)
\\ \nonumber &\leq& \psi(8c^2(1+c)).
\eeq
Therefore, we obtain (\ref{s-2.1}).

Next, to conclude the proof, we choose an arc $\alpha_0$ connecting $x$ and $y_0$ in $G$ such that $\ell_k(\alpha_0)\leq 2k(x,y_0)$. Moreover, let $y\in\alpha_0$ be the first point such that $d_G(y)=2d_G(x)$. Then Lemma \ref{s-2} follows from the choice of $\alpha=\alpha_0\cup\gamma[y_0,x_0]$.
\epf

\begin{lem}\label{s-3}  Assume that $G\varsubsetneq E$ is a $\psi$-natural and diameter $c$-uniform domain, then $G$ is $c_1$-John for some $c_1=c_1(\psi,c)$.
\end{lem}
\bpf Fix $x, y\in G$. Because $G$ is diameter $c$-uniform, there is a diameter $c$-uniform arc $\gamma$ in $G$ connecting $x$ and $y$. Choose $x_0\in\gamma$ such that $\diam \;\gamma[x,x_0]=\diam \;\gamma[x_0,y].$ Then we have
\be\label{s-3.1} \diam \;\gamma[x,x_0]=\diam \;\gamma[x_0,y]\leq cd_G(x_0).\ee
In order to show that $G$ is $c_1$-John, by symmetry, we only need to show that there is a constant $c_1>0$ and a curve $\sigma$ joining $x$ and $x_0$ such that
\be\label{s-3.2}\ell(\sigma[x,z])\leq c_1d_G(z),\ee
for all $z\in \sigma$.

We divide the proof of (\ref{s-3.2}) into two cases. Let $a_1=2\psi(8c^2(1+c))$, $c_0=2e^{a_1}+c$ and $a_2=2\psi(2cc_0(1+c))$.

\textbf{Case $1$.} Let $d_G(x)\geq d_G(x_0)/2$. By (\ref{s-3.1}) we have
$$|x-x_0|\leq cd_G(x_0).$$
Moreover, from Lemma \ref{s-1} we see that there is an arc $\sigma$ joining $x$ and $x_0$ such that
\be\label{s-3.3}\ell_k(\sigma)\leq 2\psi(2c^2(1+c))\leq a_1.\ee
Then for all $z\in\sigma$, we obtain by (\ref{eq(0000)}) and (\ref{s-3.3}) that
$$\Big|\log\frac{d_G(z)}{d_G(x_0)}\Big|\leq k(z,x_0)\leq \ell_k(\sigma)\leq a_1,$$
and therefore
$$e^{-a_1}d_G(x_0)\leq d_G(z) \leq e^{a_1}d_G(x_0).$$
This yields
\be\label{s-3.3-1} \ell(\sigma)\leq e^{\ell_k(\sigma)}d_G(x_0)\leq e^{a_1}d_G(z),\ee
as desired.

\textbf{Case $2$.} Let $d_G(x)<d_G(x_0)/2$. Let $n\geq 1$ be the unique integer such that
$$2^{n-1}d_G(x)<d_G(x_0)\leq 2^n d_G(x).$$
We  construct a sequence of points $x=x_1,x_2,\ldots,x_n,x_{n+1}=x_0$ and curves $\beta_i$ as follows.

First, we see from Lemma \ref{s-2} that there is an arc $\alpha_1$ joining $x_1$ and $x_0$ with $\ell_k(\alpha_1[x_1,x_2])\leq a_1$, where $x_2$ is the first point of $\alpha_1$ (when traversing $\alpha_1$ from $x_1$ towards $x_0$) with $d_G(x_2)=2d_G(x_1)$. If $d_G(x_2)\geq d_G(x_0)/2$, again we connect $x_2$ to $x_0$ by an arc $\alpha_2$ such that Lemma \ref{s-2} holds for this arc. Then we stop with $n=2, \beta_2=\alpha_2$ and $x_3=x_0$. Otherwise we continue the process by letting $\beta_i=\alpha_i[x_i,x_{i+1}]$ where $x_{i+1}$ is the first point of $\alpha_i$ with $d_G(x_{i+1})=2d_G(x_i)$, and where $\alpha_i$ is the arc joining $x_i$ to $x_0$ such that Lemma \ref{s-2} holds. Since $d_G(x_i)=2^{i-1}d_G(x_1)$, we find that $d_G(x_i)\geq d_G(x_0)/2$ as soon as $i\geq \log_2\frac{d_G(x_0)}{d_G(x_1)}$. Thus the above process stops with $i=n$, $\beta_n=\alpha_n$ and $x_{n+1}=x_0$.

Furthermore, we shall verify that the curve $\beta=\bigcup_{i=1}^n \beta_i$ satisfies the desired property. To this end, note that by Lemma \ref{s-2} we have
\be\label{s-3.4}\ell_k(\beta_i)\leq a_1,\ee
for $1\leq i\leq n-1$, By using a similar argument as in (\ref{s-3.3-1}), this yields
\be\label{s-3.5}\ell(\beta_i)\leq e^{a_1}d_G(x_i).\ee
For $i=n$, we claim that
\be\label{s-3.6}\ell_k(\beta_n)\leq a_2.\ee
Indeed, because $|x-x_0|\leq cd_G(x_0)$, by (\ref{s-3.5}), we have
\beq\nonumber
|x_n-x_0| &\leq& |x_n-x_{n-1}|+\ldots+|x_2-x_1|+|x_1-x_0|
\\ \nonumber &\leq& \sum_{i=1}^{n-1}\ell(\beta_i)+|x-x_0|
\\ \nonumber &\leq& \sum_{i=1}^{n-1} e^{a_1}d_G(x_i)+|x-x_0|
\\ \nonumber &\leq& (2e^{a_1}+c)d_G(x_0)=c_0d_G(x_0).
\eeq
This together with Lemma \ref{s-1} implies (\ref{s-3.6}).

Then from (\ref{s-3.6}) and a similar argument with (\ref{s-3.3-1}) it follows that
\be\label{s-3.7} \ell(\beta_n)\leq e^{a_2}d_G(x_n).\ee

Without loss of generality we may assume $a_2\geq a_1$. Fix $z\in\beta$, then there is some index $1\leq j\leq n$ such that $z\in\beta_j$. Moreover, we obtain by (\ref{eq(0000)}) and (\ref{s-3.4})
$$\Big|\log\frac{d_G(a_j)}{d_G(z)}\Big|\leq k(z,a_j)\leq \ell_k(\beta_j)\leq a_2$$
and thus
\be\label{s-3.8} d_G(x_j) \leq e^{a_2}d_G(z).\ee
Therefore, we see from (\ref{s-3.5}), (\ref{s-3.7}) and (\ref{s-3.8}) that
$$\ell(\beta[x,z])\leq \sum_{i=1}^j\ell(\beta_j)\leq e^{a_2} \sum_{i=1}^j d_G(x_j)\leq 2e^{a_2}d_G(x_j)\leq 2e^{2a_2}d_G(z),$$
which shows (\ref{s-3.2}). Hence the proof of Lemma \ref{s-3} is complete.
\epf

Finally, we are in a position to show the sufficiency part of Theorem \ref{thm-main2}. Note that by Lemma \ref{s-3}, we only need to prove the following result.

\begin{lem}\label{s-4} Assume that $G\varsubsetneq E$ is a $\psi$-natural, $c_1$-John and diameter $c$-uniform domain, then $G$ is $b$-uniform for some $b=b(\psi,c_1,c)$.
\end{lem}
\bpf Fix $x,y\in G$. Let $k$ be the quasihyperbolic metric of $G$. If $|x-y|< d_G(x)$, thus the line segment $[x,y]$ connecting $x$ and $y$ in $G$ is the desired uniform arc.

In the following, we may assume that $|x-y|\geq d_G(x)$. Since $G$ is $c_1$-John, there is a $c_1$-cone arc $\alpha$ connecting $x$ and $y$ in $G$. Choose two points $x'$ and $y'$ in $\alpha$ such that
\be\label{s-4.1}  \frac{1}{2}|x-y|=\ell(\alpha[x,x'])=\ell(\alpha[y,y'])\leq c_1\;\min\{d_G(x'),d_G(y')\}.\ee
Then by (\ref{s-4.1}) we have
\be\label{s-4.2} |x'-y'|\leq |x'-x|+|x-y|+|y-y'|\leq 2|x-y|\leq 4c_1\;\min\{d_G(x'),d_G(y')\}.\ee

Moreover, we may join $x'$ and $y'$ by a diameter $c$-uniform arc $\beta$ in $G$ because $G$ is a diameter $c$-uniform domain. It follows from the uniformity of $\beta$ and (\ref{s-4.2}) that
\be\label{s-4.3} \diam\;\beta \leq c|x'-y'|\leq 4c_1c \;\min\{d_G(x'),d_G(y')\}\ee
and
\be\label{s-4.4} \;\min\{d_G(x'),d_G(y')\}\leq d_G(z)+\min\{|x'-z|,|z-y'|\}\leq (1+c)d_G(z),\ee
for all $z\in \beta$. Because $G$ is $\psi$-natural, we obtain by (\ref{s-4.3}) and (\ref{s-4.4}) that
\be\label{s-4.5} k(x',y')\leq \diam_k(\beta)\leq \psi\Big(\frac{\diam\;\beta}{\dist(\beta,\partial G)}\Big)\leq \psi(4c_1c(1+c))=:\lambda.\ee

It thus follows from (\ref{s-4.5}) that there is an arc $\gamma$ connecting $x'$ and $y'$ in $G$ such that
\be\label{s-4.6}  \ell_k(\gamma)\leq 2k(x',y')\leq 2\lambda.\ee

Denote $\sigma=\alpha[x,x']\cup \gamma\cup \alpha[y',y]$. We show that $\sigma$ satisfies the uniformity condition. Because $d_G(x)\leq |x-y|$, by (\ref{s-4.1}) we have
\be\label{s-4.7} d_G(x')\leq d_G(x)+|x'-x|\leq 2|x-y|.\ee
Moreover, we see from (\ref{s-4.6}) that
$$\log\Big(1+\frac{\ell(\gamma)}{d_G(x')}\Big)\leq \ell_k(\gamma)\leq 2\lambda,$$
which together with (\ref{s-4.7}) implies that
$$\ell(\gamma)\leq 2e^{2\lambda}|x-y|.$$
Therefore, we obtain by (\ref{s-4.1}) that
\be\label{s-4.8} \ell(\sigma)=\ell(\alpha[x,x'])+\ell(\gamma)+\ell(\alpha[y',y])\leq (1+2e^{2\lambda})|x-y|.\ee

It is left to show that $\sigma$ satisfies the cone condition. Fix $z\in \sigma$. If $z\in \alpha[x,x']\cup \alpha[y',y]$, then the desired cone property follows from the fact that $\alpha$ is $c_1$-cone. If $z\in \gamma$, again by (\ref{eq(0000)}) and (\ref{s-4.6}) we have
$$\Big|\log\frac{d_G(x')}{d_G(z)}\Big|\leq k(z,x')\leq \ell_k(\gamma)\leq 2\lambda,$$
and therefore $d_G(x')\leq e^{2\lambda}d_G(z).$ By (\ref{s-4.8}) and (\ref{s-4.1}), it follows that
$$\ell(\sigma)\leq  (1+2e^{2\lambda})|x-y|\leq 2c_1(1+2e^{2\lambda})d_G(x')\leq 2c_1(1+2e^{2\lambda})e^{2\lambda}d_G(z),$$
as desired. This proves Lemma \ref{s-4}.
\epf

\section{Proofs of Theorems \ref{thm-main3} and \ref{thm-main4}}\label{sec-4}
In this section, we start by proving certain auxiliary results.  We first show that the diameter uniformity of domains is preserved under quasisymmetric mappings.
\begin{lem}\label{r-4.0} Let $G$ and $G'$ be proper domains in Banach spaces $E$ and $E'$, respectively. Suppose that $G$ is diameter $c$-uniform and that a homeomorphism ${f}:{G}\to {G'}$ is $\eta$-quasisymmetric, then $G'$ is diameter $c'$-uniform with $c'=c'(c,\eta)$.
\end{lem}
\bpf
For any pair of points $x',y'$ in $G'$, denote $x=f^{-1}(x')$ and $y=f^{-1}(y')$. Because $G$ is diameter $c$-uniform, we know that there is a diameter $c$-uniform arc $\gamma$ connecting $x$ and $y$ in $G$. We show that the image arc $f(\gamma)=\gamma'$ is the desired uniform arc in $G'$. Without loss of generality, we may assume that ${f}:\overline{G}\to\overline{ G'}$ is also $\eta$-quasisymmetric by \cite[Theorem 6.12]{Vai-5}.

On the one hand, for all $u\in \gamma,$ the diameter uniformity of $\gamma$ implies that
$$|x-u|\leq c|x-y|.$$
Because $f$ is $\eta$-quasisymmetric, we see that
$$|f(x)-f(u)|\leq \eta(c)|f(x)-f(y)|,$$
and thus
\be\label{r-4.2} \diam\; \gamma'\leq 2\eta(c) |f(x)-f(y)|.\ee

On the other hand, for any $z'\in \gamma'$ with $f(z)=z'$, without loss of generality we may assume that  $$\diam\;\gamma[x,z]\leq \diam\;\gamma[y,z].$$
Because $\gamma$ is diameter $c$-uniform, we have for all $v\in \gamma[x,z]$ and for every $w\in \partial G$,
$$|z-v|\leq c|z-w|.$$
By the quasisymmetry proeprty of $f$, we obtain
$$|f(z)-f(v)|\leq \eta(c)|f(z)-f(w)|.$$
This yields
\be\label{r-4.3} \diam \;\gamma'[x',z']\leq 2\eta(c) d_G(z').\ee

Therefore, from (\ref{r-4.2}) and (\ref{r-4.3}) we see that $\gamma'$ is as required.
Hence the proof of Lemma \ref{r-4.0} is complete.
\epf

Next, we prove that  the $\delta$-uniformity condition is preserved under quasim\"obius maps.

\begin{lem}\label{r-5.0} Let $G$ and $G'$ be proper domains in Banach spaces $E$ and $E'$, respectively. Let $0<\delta<1$. Suppose that $G$ is $\delta$-uniform and that a homeomorphism ${f}:{G}\to {G'}$ is $\eta$-quasim\"obius, then $G'$ is $\delta'$-uniform with $\delta'=\delta'(\delta,\eta)\in (0,1)$.\end{lem}
\bpf We want to reduce the situation to the quasisymmetric case by using auxiliary inversions. To this end, consider the one-point extensions $\dot{E}=E\cup\{\infty\}$ and $\dot{E'}=E'\cup\{\infty\}$. For all $x\in \dot{E}$, we define the inversion $u:\dot{E}\to \dot{E}$ as
$$u(x)=\frac{x}{|x|^2},$$
with $u(0)=\infty$ and $u(\infty)=0$. Similarly, for every $x'\in \dot{E'}$, the inversion $u':\dot{E'}\to \dot{E'}$ is defined by
$$u'(x')=\frac{x'}{|x'|^2},$$
with $u'(0)=\infty$ and $u'(\infty)=0$. It follows from \cite[Theorem 6.22]{Vai-5} that $u$ and $u'$ are both $\theta$-quasim\"obius with $\theta(t)=81t$. Since $G$ is $\delta$-uniform, we see that $u(G)$ is $\delta_1$-uniform with $\delta_1=\delta/81$. Similarly, $G'$ is $\delta'$-uniform if and only if $u'(G')$ is $\delta_1'$-uniform with the constants $\delta'$ and $\delta_1'$ depending only on each other.

By \cite[Theorem 6.24]{Vai-5} we may assume that $f:\overline{G}\to \overline{G'}$ is $\eta$-quasim\"obius as well. By auxiliary translations, we may also assume that $0\in\partial G$ and $f(0)=0\in \partial G'$, or $f(0)=\infty\in \partial G'$. Now we consider two possibilities. If $f(0)=0$, we define
$$g:=u'\circ f\circ u^{-1}: u(G)\to u'(G').$$
If $f(0)=\infty$, we define
$$g:= f\circ u^{-1}: u(G)\to G'.$$
Thus, we know that $g$ is $\eta_1$-quasim\"obius with $\eta_1(t)=81\eta(81t)$. Because in both cases $g(\infty)=\infty$, we see from \cite[Theorem 3.20]{Vai-0} that $g$ is $\eta_1$-quasisymmetric.

Furthermore, because $u(G)$ is $\delta_1$-uniform, by Lemma \ref{r-3.0.1} we find that $u(G)$ is diameter $c$-uniform with $c=c(\delta_1)$. It thus follows from Lemma \ref{r-4.0} that $g\circ u(G)$ is diameter $c'$-uniform with $c'=c'(c,\eta)$. Therefore by using Lemma \ref{r-3.0}, we obtain that $G'$ is $\delta'$-uniform with $\delta'=\delta'(\delta,\eta)$, as desired.

Hence Lemma \ref{r-5.0} is proved.
\epf

\begin{lem}\label{r-6.0} Let $G$ and $G'$ be proper domains in Banach spaces $E$ and $E'$, respectively.  Suppose that $G$ is  diameter $c$-uniform  and that a homeomorphism ${f}:{G}\to {G'}$ is $\eta$-quasim\"obius. Then $G'$ is diameter $c'$-uniform with $c'=c'(c,\eta)$.
\end{lem}
\bpf This lemma follows immediately from Lemmas \ref{r-3.0}, \ref{r-3.0.1} and \ref{r-5.0}.
\epf

\begin{lem}\label{r-7.0} Let $G$ and $G'$ be proper domains in Banach spaces $E$ and $E'$, respectively. Suppose that $G$ has the min-max property and that a homeomorphism ${f}:\overline{G}\to \overline{G'}$ is $\eta$-quasisymmetric relative to  $\partial G$ and maps $G$ onto $G'$. Then $G'$ also has the min-max property.
\end{lem}
\bpf Suppose that there exists a family of curves $\Gamma$ in $G$ and a constant $c\geq 1$ such that any pair of points in $G$ can be joined by a curve $\gamma\in \Gamma$ and so that (\ref{r-1}) holds for each ordered triplet of points $x_1,x,x_2\in \gamma$ and for all $y\in \partial G$.

Because $f$ is $\eta$-quasisymmetric relative to  $\partial G$, it is not difficult to see that any pair of points in $G'$ can be joined by a curve $\gamma'\in \Gamma'=f(\Gamma)$ such that (\ref{r-1}) holds with the constant $\eta(c)$, for each ordered triplet of points $x_1',x',x_2'\in \gamma'$ and for all $y'\in \partial G'$. This completes the proof of Lemma \ref{r-7.0}.
\epf

Now, we are ready to complete the proofs of Theorems \ref{thm-main3} and \ref{thm-main4}.

\textbf{Proof of Theorem \ref{thm-main3}.} We know from Theorem \ref{thm-main2} that the second assertion of Theorem \ref{thm-main3} follows from the first statement. Thus we only need to show that $G'$ is diameter uniform. By using a similar argument as Lemma \ref{r-5.0}, we reduce the situation to the case of relative quasisymmetry by using auxiliary inversions. For completeness, we show the details.

Consider the one-point extensions $\dot{E}=E\cup\{\infty\}$ and $\dot{E'}=E'\cup\{\infty\}$. By auxiliary translations, there is no loss of generality in assuming that $0\in\partial G$ and $f(0)=0\in \partial G'$, or $f(0)=\infty\in \partial G'$. Now we consider two possibilities. If $f(0)=0$, we define
$$g:=u'\circ f\circ u^{-1}: u(G)\to u'(G').$$
If $f(0)=\infty$, in this case we define
$$g:= f\circ u^{-1}: u(G)\to G'.$$

By \cite[Theorem 6.22]{Vai-5}, we know that $g$ is $\eta_1$-quasim\"obius relative to  $\partial(u(G))$ with $\eta_1(t)=81\eta(81t)$. Because in both cases $g(\infty)=\infty$, we see that $g$ is $\eta_1$-quasisymmetric relative to  $\partial(u(G))$. Moreover, because $G$ is $c$-uniform, we see from Theorem \Ref{T-2} that $u(G)$ is $c_1$-uniform with $c_1=c_1(c)$. Thus by Lemma \ref{r-1.0}, we find that $u(G)$ satisfies the min-max property.  Furthermore, it follows from Lemma \ref{r-7.0} that $g\circ u(G)$ also satisfies the min-max property. Therefore, we may use Lemma \ref{r-2.0} to conclude that $g\circ u(G)$ is diameter $c_1'$-uniform with $c_1'$ depending only on $c$ and $\eta$.

Hence by Lemma \ref{r-6.0}, we get that $G'$ is diameter $c'$-uniform with $c'$ depending only on $c$ and $\eta$. This completes the proof of Theorem \ref{thm-main3}.
\qed

\textbf{Proof of Theorem \ref{thm-main4}.} We first prove $(1)$. Let $A'$ be a nonempty connected set in $G'$ with $r_{G'}(A')=t$ and $A=f^{-1}(A')\subseteq G$. Because $f$ is $\eta$-quasim\"obius relative to  $\partial G$, we find that the inverse map $f^{-1}:\overline{G}\to \overline{G'}$ is $\eta'$-quasim\"obius relative to  $\partial G'$ with $\eta'(t)=\eta^{-1}(t^{-1})^{-1}$ (cf. \cite{Vai-5}). Then, by Lemma \Ref{z-2}, we see that there is an increasing function $\mu:[0,\infty)\to [0,\infty)$ depending only on $\eta$ such that
\be\label{z-2.1}  r_G(A)\leq \mu(t).\ee
Moreover, because $f$ is $(M,C)$-CQH and $G$ is $\psi$-natural, we obtain by (\ref{z-2.1}) that
$$k_{G'}(A')\leq M\diam_k(A)+C\leq M\psi(r_G(A))+C\leq M\psi\circ\mu(t)+C.$$
This implies that $G'$ is $\psi'$-natural by taking $\psi'(t)=M\psi\circ\mu(t)+C$, as desired.

Next, we verify $(2)$. Assume that $G$ is $c$-uniform with $c\geq 1$. By Lemma \Ref{lem-3}, we see that $G$ is $\psi$-natural with $\psi$ depending only on $c$. Moreover, the statement $(1)$ shows that $G'$ is $\psi'$-natural. On the other hand, by Theorem \ref{thm-main3} we obtain that $G'$ is $c'$-uniform. Hence the proof of Theorem \ref{thm-main4} is complete.
\qed

\section{Proofs of Corollaries \ref{c-1}, \ref{c-2} and \ref{c-3}}\label{sec-5}
\textbf{Proof of Corollary \ref{c-1}.} Corollary \ref{c-1} follows immediately from Lemmas \ref{r-5.0} and \ref{r-6.0}.\qed

\textbf{Proof of Corollary \ref{c-2}. } Because the $\psi$-uniformity of domains implies the $\psi$-naturality, the sufficiency in Corollary \ref{c-2} follows from Theorem \ref{thm-main2}. The necessity is proved by \cite[Theorem 6.16]{Vai-2} because a $c$-uniform domain is clearly diameter $c$-uniform. The second assertion follows from the fact that a convex domain in Banach space is $\varphi$-uniform with $\varphi(t)=t$.
\qed

\textbf{Proof of Corollary \ref{c-3}.}
By Theorem \ref{thm-main4}, we see that $G$ is $\psi_1$-natural. Because $G$ is diameter $c$-uniform, we know from Theorem \ref{thm-main2} that $G$ is $c_0$-uniform. Again, it follows from Theorem \ref{thm-main4} that $G'$ is $c_1$-uniform. This proves Corollary \ref{c-3}.\qed



\end{document}